\documentclass[12pt]{amsart}
\usepackage{amssymb, amstext, amscd, amsmath}
\usepackage{fullpage}

\theoremstyle{plain}
\newtheorem{thm}{Theorem}[section]
\newtheorem{cor}[thm]{Corollary}
\newtheorem{prop}[thm]{Proposition}
\newtheorem{lem}[thm]{Lemma}

\theoremstyle{definition}
\newtheorem{defn}[thm]{Definition}

\def\proofbegin{\hspace{-\parindent}{\bf {}Proof}. }
\newcommand{\proofend}{\hspace*{\fill}$\blacksquare$\hbox{\quad}
\medskip
}

\theoremstyle{remark}

\begin{document}

\title{Generalized notions of character amenability}

\author[L. Y. Shi]{Luo Yi Shi}

\address{Department of Mathematics\\Tianjin Polytechnic University\\Tianjin 300160\\P.R. CHINA}

\email{sluoyi@yahoo.cn}

\author[Y. J. Wu]{YU Jing Wu}

\address{Tianjin Vocational Institute \\Tianjin 300160\\P.R. CHINA}

\email{wuyujing111@yahoo.cn}

\author[Y.Q. Ji]{You Qing Ji}

\address{Department of Mathematics\\Jilin University\\Changchun 130012\\P.R. CHINA}

\email{jiyq@jlu.edu.cn}

\thanks{}

\date{}



\begin{abstract}

In this paper  the concepts of character contractibility,
approximate character
 amenability (contractibility) and uniform approximate
character
 amenability (contractibility) are  introduced. We
 are concerned with the relations among the generalized concepts of character
 amenability for Banach algebra. We prove that  approximate
 character
amenability and approximate character contractibility are the same
properties, as are uniform approximate character amenability and
character amenability, as are uniform approximate character
contractibility and character contractibility. For commutative
Banach algebra, we prove that character contractibility and
contractibility are the same properties. Moreover,
 general theory for  those concepts is developed.
\end{abstract}

\maketitle

\section{Introduction}

The concept of amenability for Banach algebras was first introduced
by B. E. Johnson in \cite{BE1972}. Suppose that $A$ is a Banach
algebra and that $E$ is a Banach $A$-bimodule, then $E^*$, the dual
of $E$, has a natural Banach $A$-bimodule structure defined by
 $$ (a \cdot f)(x)=f(x \cdot a), (f \cdot a)(x)=f(a \cdot x), a\in A, x\in E,
f\in E^*.$$
 Such a Banach $A$-bimodule $E^*$ is called a {\it dual
$A$-bimodule}.
 A {\it derivation} $D : A\rightarrow E^*$ is a continuous linear map such that $D(ab) =
a \cdot D(b) + D(a)\cdot b$  for all $a, b \in A$. Given $f\in E^*$,
the {\it inner derivation} $\delta_f : A\rightarrow E^*$, is defined
by $\delta_f (a) = a\cdot f- f \cdot a$. According to Johnson's
original definition, a Banach algebra $A$ is {\it amenable} if every
 derivation from $A$ into the dual  $A$-bimodule
$E^*$ is inner for all Banach $A$-bimodules $E$.
 As a complement to this notion, a Banach
algebra $A$ is {\it contractible} if every derivation from $A$ into
every Banach $A$-bimodule is inner \cite{P1989,A1993}.

Ever since its introduction, the concept of amenability has occupied
an important place in the research of Banach algebras, operator
algebras and harmonic analysis. For example,  an early result of
Johnson \cite{V1774} shows that the amenability of the group algebra
$L^1(G)$, for G a locally compact group, is equivalent to the
amenability of the underlying group G. Results of Connes and
Haagerup show that a $C^*$-algebra is amenable if and only if it is
nuclear \cite{V1774}. However it has been realized that  amenability
is  essentially a  finiteness condition, and in many instances is
too restrictive. As for contractibility, it is even conjectured in
\cite{A1993} that a contractible Banach algebra must be finite
dimensional (see also \cite{Y2000}). For this reason by relaxing
some of the constrains in the definition of amenability new concepts
have been introduced. The most notable are the concepts of Connes
amenability \cite{A1991, B1972}, weak amenability \cite{W1987,H2000}
and character amenability \cite{EAJ2008, M2008}. More recently, F.
Ghahramani and R. J. Loy have introduced and studied the concepts of
approximate amenability (contractibility) and uniform approximate
amenability (contractibility) for Banach algebras
\cite{G2004,G2008}. In this paper we introduce the generalized
concepts of character amenability (see Definition 1.2-1.7). We are
concerned with the relations among those concepts and shall
 develop general theory for them.

\begin{defn}
A derivation $D : A\rightarrow E$ is {\it approximately inner}, if
there exists a net $\{\xi_i\}\subset E$ such that
$D(a)=\lim\limits_i(a\cdot\xi_i-\xi_i\cdot a)$ for all $a\in A$, the
limit being in norm.
\end{defn}

Note that $\{\xi_i\}$ in the above is not necessarily  bounded. The
stronger assumption, that $D$ is in the uniform closure of the inner
derivations, has been well studied in the $C^*$-algebra case with
the restriction to the single Banach $A$-bimodule $E=A$ (see
\cite{R1978,R1967}). The case of semigroup algebras is considered in
\cite{S1998} for the Banach $A$-bimodule $A^*$.

 Let $A$ be a Banach algebra and
$\sigma(A)$ be the set of all non-zero multiplicative linear
functionals on $A$. If $\varphi \in \sigma(A)\cup\{0\}$ and $E$ is
an arbitrary Banach space, then $E$ can be viewed as a Banach left
or right $A$-module by the following actions. For $a\in A, x\in E$:
$$ a\cdot x = \varphi(a)x,\eqno(2.1)$$
$$ x\cdot a =\varphi(a)x.\eqno(2.2)$$
If the left action of $A$ on $E$ is given by (2.1), then it is
easily verified that the right action of $A$ on the dual $A$-module
$E^*$ is given by $f\cdot a = \varphi(a) f$ for all $f\in E^*, a\in
A$. Throughout, by a {\it $(\varphi, A)$-bimodule} $E$, we mean that
$E$ is a Banach $A$-bimodule for which the left module action is
given by (2.1). {\it $(A, \varphi)$-bimodule} is defined similarly
by (2.2).

\begin{defn}
 Let $A$ be a Banach algebra and $\varphi\in
\sigma(A)\cup\{0\}$.  $A$ is {\it approximately $\varphi$-amenable},
 if every derivation $D$ from $A$ into the dual
$A$-bimodule $E^*$ is approximately inner for all  $(\varphi,
A)$-bimodules $E$.
\end{defn}

\begin{defn}
  A  Banach algebra $A$ is
{\it approximately right character amenable}, if for every
$\varphi\in \sigma(A)\cup\{0\}$ and every  $(\varphi, A)$-bimodule
$E$, every derivation $D$ from $A$ into the dual $A$-bimodule $E^*$
is approximately inner.
\end{defn}

{\it Approximately left character amenability} is defined similarly
by considering  $(A, \varphi)$-bimodules $E$.

$A$ is {\it approximately character amenable} if it is both
approximately left and right character amenable.

\begin{defn}
 Let $A$ be a Banach algebra and $\varphi\in
\sigma(A)\cup\{0\}$.  $A$ is {\it $\varphi$-contractible},
 if  every derivation $D : A\rightarrow E$ is inner, for all
  $(A, \varphi)$-bimodules
$E$.
\end{defn}

\begin{defn}
A  Banach algebra $A$ is {\it left character contractible}, if  for
every $\varphi\in \sigma(A)\cup\{0\}$ and every  $(\varphi,
A)$-bimodule $E$, every derivation $D : A\rightarrow E$ is inner.
\end{defn}

{\it Right character contractibility} is defined similarly by
considering  $(A, \varphi)$-bimodules $E$.

 $A$ is {\it character contractible}
if it is both left and right character contractible.

\begin{defn}
 Let $A$ be a Banach algebra and $\varphi\in
\sigma(A)\cup\{0\}$.  $A$ is {\it approximately
$\varphi$-contractible},
 if  every derivation $D: A\rightarrow E$ is
approximately inner for all  $(A, \varphi)$-bimodules $E$.
\end{defn}

\begin{defn}
 A  Banach algebra $A$ is {\it approximately
left character contractible}, if  for every  $\varphi\in
\sigma(A)\cup\{0\}$ and every  $(\varphi, A)$-bimodule $E$, every
 derivation $D : A\rightarrow E$ is approximately inner.
\end{defn}

{\it Approximately right character contractible} is defined
similarly by considering  $(A, \varphi)$-bimodules $E$.

 $A$ is {\it approximately character contractible}
if it is both approximately  left  and right character contractible.

Any statement about approximate left character amenability
(contractibility) turns into an analogous statement about
approximate right character amenability (contractibility,
respectively) by simply replacing $A$ by its opposite algebra.
Approximate right character amenability (contractibility) of $A$ is
equivalent to approximate $\varphi$-amenability
($\varphi$-contractibility) for all $\varphi\in\sigma(A)$ together
with approximate $0$-amenability ($0$-contractibility,
respectively).

The qualifier {\it uniform} on the above definitions will indicate
that the convergence of the net is uniform over the unit ball of
$A$. Similarly {\it $w^*$} will indicate that  convergence is in the
appropriate $w^*$-topology. Let {\bf $CC$} denote character
contractibility, {\bf $UACC$} denote uniform approximate character
contractibility, {\bf $ACC$} denote approximate character
contractibility, {\bf $CA$} denote character amenability, {\bf
$UACA$} denote uniform approximate character amenability and {\bf
$ACA$} denote approximate character amenability. Clearly the
relations among the various character amenability are as follows:
\begin{center}
$Contractibility \Rightarrow CC \Rightarrow  UACC \Rightarrow ACC$

$\Downarrow \ \ \ \ \not\Uparrow \ \ \ \ \ \ \ \ \ \ \ \ \ \
\Downarrow \ \ \ \ \ \ \ \ \ \ \ \Downarrow \ \ \ \ \ \ \ \
\Downarrow$

$Amenability \Rightarrow CA \Rightarrow UACA \Rightarrow ACA$
\end{center}
In this paper, we prove that $ACA\Leftrightarrow ACC$ (Theorem
\ref{T:Equi 1}), $UACA\Leftrightarrow CA$ (Theorem \ref{T: Uaa-a}),
$UACC\Leftrightarrow CC$ (Theorem \ref{T: Uac-c}). Moreover,  in
Section 6, Example 1 shows that $ACA\nRightarrow UACA$.  Example 2
shows that $CA\nRightarrow Amenability$  and $CA\nRightarrow CC$.
Thus $ACC\nRightarrow UACC$.  For commutative Banach algebra, we
obtain that $Contractibility\Leftrightarrow CC$ (Theorem \ref{T:
c-cc}). For non-commutative Banach algebra, Example 3 shows that
$CC\nRightarrow Amenability$.

In this paper, the second dual $A^{**}$ of a Banach algebra $A$ will
always be equipped with the first Arens product  \cite{J1979} which
is defined as follows. For $a, b \in A, f \in A^*$ and $m,n \in
A^{**}$, the elements $f\cdot a$ and $m\cdot f$ of $A^*$ and $mn\in
A^{**}$ are defined by $$(f\cdot a)(b)=f (ab), (m\cdot
f)(b)=m(f\cdot b), mn(f) = m(n\cdot f ),$$ respectively. With this
multiplication, $A^{**}$ is a Banach algebra and $A$ is a subalgebra
of $A^{**}$.
 Moreover, for all
$m,n\in A^{**}$ and $\varphi\in\sigma(A)$,
$(mn)(\varphi)=m(\varphi)n(\varphi)$. Consequently, each
$\varphi\in\sigma(A)$ extends uniquely to some element
$\varphi^{**}$ of $\sigma(A^{**})$. The kernel of $\varphi^{**}$,
$\ker(\varphi^{**})$, contains $\ker\varphi$ in the same sense that
$A^{**}$ naturally contains $A$. Since each of these ideals has
codimension 1, the theory of second dual shows that $\ker\varphi$ is
$w^*$-dense in $\ker(\varphi^{**})$ and that $\ker\varphi^{**} =
(\ker\varphi)^{**}$. For further details  the reader is referred to
\cite{J1979}.

The organization of the paper is as follows. In Section 2 we
characterize (uniform) approximate character amenability in three
different ways. In Section 3 we are concerned with  hereditary
properties of (uniform) approximate character amenability.

In Section 4 we characterize character contractibility and (uniform)
approximate character contractibility in two different ways.

Section 5 is devoted to the  relations among generalized notions of
character-amenability. We prove that  approximate
 character
amenability and approximate character contractibility are the same
properties, as are uniform approximate character amenability
(contractibility) and character amenability (contractibility,
respectively). For commutative Banach algebra, character
contractibility and contractibility are the same properties.

Section 6 gives three examples. The first example shows that there
exists a Banach algebra which is approximately character amenable
but not uniformly approximately character amenable. The second
example shows that there exists  a Banach algebra which is character
amenable but neither amenable nor character contractible. The last
example shows that there exists a Banach algebra which is character
contractible but not amenable.

\vskip1cm
\section{Characterization of approximately character amenability}

In this section, we first characterize (uniform) approximate
$\varphi$-amenability in three
 different ways and then characterize (uniform) approximate character amenability in these different ways.

Suppose that $A$ is a Banach algebra, we make $E=A$ into a Banach
$A$-bimodule as usual by $a\cdot b=ab, b\cdot a=ba,$ for all $a\in
A, b\in E$. Then $A^*, A^{**}$ are dual $A$-bimodules and the module
actions are given by $(a\cdot f)(b)=f(ba), (f\cdot a)(b)=f(ab);
(a\cdot m)(f)=m(f\cdot a), (m\cdot a)(f)=m(a\cdot f)$ for all $a,
b\in A, f\in A^*, m\in A^{**}$. If we take $A^{**}$ with the first
Arens product, then $am=a\cdot m, ma=m\cdot a $ for all $a\in A,
m\in A^{**}$. Let $\varphi\in \sigma(A)$, a net $\{m_\alpha\}\subset
A^{**}$ is called an {\it approximate $\varphi$-mean}, if
$m_\alpha(\varphi)=1$ and $||a\cdot
m_\alpha-\varphi(a)m_\alpha||\rightarrow 0$ for all $a\in A$. The
following proposition characterizes approximate $\varphi$-amenable
in terms of approximately $\varphi$-mean. The corresponding result
characterizing $\varphi$-amenability of a Banach algebra was
obtained in \cite[ Theorem 1.1]{E2008}.

\begin{prop}\label{P:P-A1}
 Let $A$ be a Banach algebra and
$\varphi\in\sigma(A)$. Then the following are equivalent:

$(i)$ $A$ is (uniformly) approximately $\varphi$-amenable;

$(ii)$ There exists a net $\{ m_\alpha\}\subset A^{**}$ such that
$m_\alpha(\varphi)=1$ and $||a\cdot
m_\alpha-\varphi(a)m_\alpha||\rightarrow 0$, for all $a\in A$
(uniformly on the unit ball of $A$, respectively);

$(iii)$ There exists a net $\{ m_\alpha\}\subset A^{**}$ such that
$m_\alpha(\varphi)\rightarrow1$ and $||a\cdot
m_\alpha-\varphi(a)m_\alpha||\rightarrow 0$, for all $a\in A$
(uniformly on the unit ball of $A$, respectively);

$(iv)$  Give $(\ker\varphi)^{**}$ a dual $A$-bimodule structure by
taking the right action to be $m\cdot a =\varphi(a)m$ for $m\in
A^{**}$ and taking the left action to be the natural one. Then any
continuous derivation $D : A\rightarrow(\ker\varphi)^{**}$is
(uniformly, respectively) approximately inner.
\end{prop}

\proofbegin $(i)\Rightarrow(iv)$ and $(ii)\Rightarrow(iii)$ is
clear. Therefore, in order to establish the proposition it suffices
to show the implications $(i)\Rightarrow(ii)$, $(iii)\Rightarrow(i)$
and $(iv)\Rightarrow(ii)$.

$(i)\Rightarrow(ii)$ We  first define an action of $A$ on $E=A^*$ by
$$a*f=\varphi(a)f, f* a=f\cdot a,~ a\in A, f\in E,$$
then $E$ is $(\varphi, A)$-bimodule and $A^{**}$ is a dual
$A$-bimodule and module actions are defined by $a* m=a\cdot m, m*a
=\varphi(a)m$, for all $a\in A, m\in A^{**}$.

We know that $\varphi\in A^*$, and $a\cdot \varphi=\varphi\cdot
a=\varphi(a)\varphi$. Therefore $\mathbb{C}\varphi=\{\lambda\varphi,
\lambda\in \mathbb{C}\}$ is a closed submodule of $A^*$ and
$A^*/\mathbb{C}\varphi$ is a $(\varphi, A)$-bimodule for which the
module actions are given by $a\cdot [f]=\varphi(a)[f], [f]\cdot
a=[f*a]$, for all $a\in A, [f]\in A^*/\mathbb{C}\varphi$.

Choose any $m\in A^{**}$ with $m(\varphi)=1$, and define a
derivation $D:A\rightarrow A^{**}$ by $D(a)=a\cdot m-\varphi(a)m$,
then $D(a)\in \{n\in A^{**},
n(\varphi)=0\}=\{\mathbb{C}\varphi\}^\bot\cong
(A^*/\mathbb{C}\varphi)^*$. Since $A$ is approximately
$\varphi$-amenable, it follows that there exists a  net
$\{n_\alpha\}\subset\{\mathbb{C}\varphi\}^\bot$ such that
$D(a)=\lim\limits_\alpha(a\cdot n_\alpha-\varphi(a)n_\alpha)$. Set
$m_\alpha=m-n_\alpha$, then $m_\alpha(\varphi)=1$ and $||a\cdot
m_\alpha-\varphi(a)m_\alpha||\rightarrow 0$ for all $a\in A$.

$(iii)\Rightarrow(i)$  Let $\{ m_\alpha\}\subset A^{**}$ such that
$m_\alpha(\varphi)\rightarrow1$ and $||a\cdot
m_\alpha-\varphi(a)m_\alpha||\rightarrow 0$ for all $a\in A$.   Let
$E$ be a $(\varphi, A)$-bimodule.  Also, let $D: A\rightarrow E^*$
be a continuous derivation, and let $D^{'}=D^*|_{E}: E\rightarrow
A^*$ and $g_\alpha=(D^ {'})^{*} (m_\alpha)\in E^*$. Then, for all
$a,b \in A$ and $x\in E$,
\begin{eqnarray*}
D^{'}(x\cdot a)(b)&=&D(b)(x\cdot a)\\
&=&D(ab)(x)-D(a)(b\cdot x)\\
&=&(D^{'}(x)\cdot a)(b)-D(a)(x)\varphi(b),
\end{eqnarray*}
and hence $D^{'}(x\cdot a)=D^{'}(x)\cdot a-D(a)(x)\varphi.$ This
implies that
\begin{eqnarray*}
(a\cdot g_\alpha)(x)&=&g_\alpha(x\cdot a)\\
&=&(D^ {'})^{*} (m_\alpha)(x\cdot a)\\
&=&m_\alpha(D^{'}(x\cdot a))\\
&=&m_\alpha(D^{'}(x)\cdot a)-D(a)(x)m_\alpha(\varphi)\\
&=&(a\cdot m_\alpha)(D^{'}(x))-D(a)(x)m_\alpha(\varphi).
\end{eqnarray*}

It follows that
\begin{eqnarray*}
||(a\cdot g_\alpha)(x)-\varphi(a)g_\alpha(x)+D(a)(x)||&\leq&
||\varphi(a)g_\alpha(x)-(a\cdot
m_\alpha)(D^{'}(x))||\\
&+&||D(a)(x)-D(a)(x)m_\alpha(\varphi)||.
\end{eqnarray*}

Hence,
\begin{eqnarray*}
\lim\limits_\alpha||(a\cdot
g_\alpha)(x)-\varphi(a)g_\alpha(x)+D(a)(x)|| &\leq&
\lim\limits_\alpha\{||\varphi(a)m_\alpha(D^{'}(x))-(a\cdot
m_\alpha)(D^{'}(x))||\\
\qquad &+&||\varphi(a)g_\alpha(x)-\varphi(a)m_\alpha(D^{'}(x))||\\
&+&||D(a)(x)-D(a)(x)m_\alpha(\varphi)||\}.
\end{eqnarray*}

Thus, for each $a\in A$, $D(a)=\lim\limits_\alpha
\varphi(a)g_\alpha-a\cdot g_\alpha$. Combining this with the
equation $g_\alpha\cdot a=\varphi(a)g_\alpha$, we obtain
$D(a)=\lim\limits_\alpha g_\alpha\cdot a -a\cdot g_\alpha$ for all
$a\in A$. Set $f_\alpha=-g_\alpha$, then $D(a)=\lim\limits_\alpha
a\cdot f_\alpha-f_\alpha\cdot a $ for all $a\in A$. Since $D$ was
arbitrary, it follows that $A$ is approximately $\varphi$-amenable.

$(iv)\Rightarrow(ii)$  Choose $b \in A$ with $\varphi(b) = 1$. Then
$Da = ab-\varphi(a)b$, $a \in A$, defines a derivation from A into
$(\ker\varphi)^{**}$. By $(iv)$, D is approximately inner, it
follows that there exists  a net
$\{n_\alpha\}\subset(\ker\varphi)^{**}$ such that $D(a)
=\lim\limits_\alpha a\cdot n_\alpha-\varphi(a)n_\alpha$ for all $a
\in A$.

Set $m_\alpha=b-n_\alpha$. Then $m_\alpha(\varphi)=1$ and $||a\cdot
m_\alpha-\varphi(a)m_\alpha||\rightarrow 0$ for all $a\in A$.

The proof in the case of uniform approximate $\varphi$-amenability
is similar. \proofend

\begin{prop}\label{P:Equi 1}
For a Banach algebra $A$ and $\varphi\in\sigma(A)$, the following
are equivalent:

$(i)$ There exists a net $\{ m_\alpha\}\subset A^{**}$ such that
$m_\alpha(\varphi)\rightarrow 1$ and $||a\cdot
m_\alpha-\varphi(a)m_\alpha||\rightarrow 0$ for all $a\in A$ ;

$(ii)$ There exists a net $\{ n_\beta\}\subset A$ such that $
\varphi(n_\beta)\rightarrow1$ and $||an_\beta-\varphi(a)
n_\beta||\rightarrow 0$ for all $a\in A$.
\end{prop}

\proofbegin It suffices to show that $(ii)\Rightarrow (i)$.

Suppose that $(ii)$ holds. Take $\varepsilon>0$ and finite sets
$F\subset A$, $\Phi\subset A^*$. Then there exists  $\alpha$ such
that
$$|(a\cdot m_\alpha-\varphi(a)m_\alpha)(f)|<\frac{\varepsilon}{3},
|m_\alpha(\varphi)|>1-\varepsilon,  a\in F, f\in\Phi.$$ By
Goldstine's theorem, there exists $b_\alpha\in A$ such that
$$|f(b_\alpha)-m_\alpha(f)|<\frac{\varepsilon}{3K},  f\in
\Phi\cup\Phi\cdot F\cup\{\varphi\},$$ where $K=\sup\{|\varphi(a)|,
a\in F\}$.

Thus, for any $f\in \Phi$ and $a\in F$,
\begin{eqnarray*}
|f(ab_\alpha-\varphi(a)b_\alpha)|&\leq&|f(ab_\alpha)-(a\cdot m_\alpha)(f)|\\
&+&|(a\cdot
m_\alpha)(f)-\varphi(a)m_\alpha(f)|\\
&+&|\varphi(a)m_\alpha(f)-\varphi(a)f(b_\alpha)|
\\
&\leq&|(f\cdot a)( b_\alpha)-m_\alpha(f\cdot a )|\\
&+&|(a\cdot
m_\alpha)(f)-\varphi(a)m_\alpha(f)|\\
&+&|\varphi(a)m_\alpha(f)-\varphi(a)f(b_\alpha)|\\
&<&\varepsilon.
\end{eqnarray*}

Then there exists a net $\{b_\lambda\}\subset A$ such that for every
$a\in A$, $ab_\lambda-\varphi(a)b_\lambda\rightarrow 0,
\varphi(b_\lambda)\rightarrow 1$ weakly in $A$.

Finally, for each finite set $F\subset A$, say $F = \{a_1,
a_2,\cdots, a_n \}$,
$$(a_1b_\lambda-\varphi(a_1)b_\lambda, a_2b_\lambda-\varphi(a_2)b_\lambda, \cdots ,
a_nb_\lambda-\varphi(a_n)b_\lambda, \varphi(b_\lambda))\rightarrow
(0, 0, \cdots, 0, 1)$$
 weakly in $A^{n}\oplus \mathbb{C}.$ Thus

$(0, 0, \cdots, 0, 1)\in
\overline{co}^{weak}\{(a_1b_\lambda-\varphi(a_1)b_\lambda, a_2
b_\lambda-\varphi(a_2)b_\lambda, \cdots ,
a_nb_\lambda-\varphi(a_n)b_\lambda, \varphi(b_\lambda))\}$.

The Hahn-Banach theorem now gives that for each $\varepsilon>0$,
there exists $u_{\varepsilon,F} \in  co\{b_\lambda\}$, such that
$||au_{\varepsilon,F}-\varphi(a)u_{\varepsilon,F}||<\varepsilon,
|\varphi(u_{\varepsilon,F})-1|<\varepsilon$ for all $a\in F$.
Therefore, there exists a net $\{ n_\beta\}\subset A$ such that $
\varphi(n_\beta)\rightarrow1$ and $||an_\beta-\varphi(a)
n_\beta||\rightarrow 0$ for all $a\in A$. \proofend

Let $A$ be a Banach algebra, $\{e_\alpha\}$ be a net of $A$. We call
$\{e_\alpha\}$ a {\it right approximate identity} for $A$, if
$||ae_\alpha-a||\rightarrow 0$ for all $a\in A$.  Left (two-sided)
approximate identity for $A$ is defined similarly.  We call
$\{e_\alpha\}$ a {\it bounded right (left, two-sided, respectively)
approximate identity} for $A$, if it is a bounded net. The next
proposition characterizes approximate $\varphi$-amenability in terms
of the existence of  right approximate identity for $\ker\varphi$.
The corresponding result characterizing $\varphi$-amenability of a
Banach algebra was obtained in \cite[Proposition 2.1]{E2008}.

\begin{lem}\label{L:1}
 Let $A$ be a Banach algebra and
$\varphi\in\sigma(A)$. If the ideal $I_\varphi=\ker\varphi$ has a
right approximate identity, then $A$ is approximately
$\varphi$-amenable.
\end{lem}

\proofbegin
Choose  $u_0\in A$ such that $\varphi(u_0)=1$. Then
 $a_0=u_0^2-u_0\in
I_\varphi$. Let $\{b_\alpha\}$ be a right approximate identity for
$I_\varphi$.

Set $m_\alpha=u_0-u_0b_\alpha\in A$. Then, for any $b\in I_\varphi$,
$$||b(u_0-u_0b_\alpha)||=||bu_0-bu_0b_\alpha||\rightarrow 0.$$

Furthermore,
$$||u_0(u_0-u_0b_\alpha)-(u_0-u_0b_\alpha)||
=||u_0^2-u_0^2b_\alpha-u_0+u_0b_\alpha||
=||a_0-a_0b_\alpha||\rightarrow 0.$$

It follows that, $\varphi(m_\alpha)=1$ and
$||am_\alpha-\varphi(a)m_\alpha||\rightarrow 0$, for all $a\in A$.
Thus $A$ is approximately $\varphi$-amenable by Proposition
\ref{P:P-A1}/\ref{P:Equi 1}. \proofend

\begin{lem}\label{L:2}
 Suppose that $A$ is approximately
$\varphi$-amenable for some $\varphi\in\sigma(A)$ and that $A$ has a
right approximate identity. Then $I_\varphi$ has a right approximate
identity.
\end{lem}

\proofbegin
 Choose $u_0\in A$ such that $\varphi(u_0)=1$ and
$A=\mathbb{C}u_0\oplus I_\varphi$. Let $n_\beta=\lambda_\beta
u_0+b_\beta$ be a right approximate identity for $A$, where
$b_\beta\in I_\varphi$ and $\lambda_\beta\rightarrow 1$. Since $A$
is approximately $\varphi$-amenable, it follows from Proposition
\ref{P:P-A1}/\ref{P:Equi 1} that there exists a net
$m_\alpha=\lambda_\alpha u_0+b_\alpha\in A$ such that
$||am_\alpha-\varphi(a)m_\alpha||\rightarrow 0$, where $b_\alpha\in
I_\varphi$ and $\lambda_\alpha\rightarrow 1$.

Set $e_{\alpha,\beta}=b_\beta-b_\alpha$. Notice  that
$||be_{\alpha,\beta}-b||\rightarrow 0,$ for any $b\in I_\varphi$. In
fact, for any $b\in I_\varphi$, $||bn_\beta-b||\rightarrow 0$ and
$||bm_\alpha||\rightarrow 0$. It follows that
$||bu_0+bb_\beta-b||\rightarrow 0$ and
$||bu_0+bb_\alpha||\rightarrow 0$. Thus
$||be_{\alpha,\beta}-b||\rightarrow 0$ for all $b\in I_\varphi$.
\proofend

The following result follows immediately from Lemma
\ref{L:1}/\ref{L:2} and we omit its proof.
\begin{prop}\label{P:RAIA}
 Let $A$ be a Banach algebra with a right
approximate identity and $\varphi\in\sigma(A)$. Then $A$ is
approximate $\varphi$-amenable if and only if $I_\varphi$ has a
right approximate identity.
\end{prop}

If $A$ is a Banach algebra and $A \widehat{\otimes}A$ denotes the
projective product \cite{V1774}, then the corresponding {\it
diagonal operator} is defined as
$$\triangle: A \widehat{\otimes}A\rightarrow A, a\otimes b\rightarrow ab.$$
Then $A \widehat{\otimes}A$ becomes a Banach $A$-bimodule through
$$a\cdot (b\otimes c)=ab\otimes c, (b\otimes c)\cdot a=b\otimes ca$$
for all $a, b, c\in A$. By \cite[Theorem 2.2.4]{V1774}, $A$ is
amenable if and only if there is a net $\{m_\alpha\}\subset A
\widehat{\otimes}A$ such that $a\cdot m_\alpha-m_\alpha\cdot
a\rightarrow 0, a\triangle m_\alpha\rightarrow a$ for all $a\in A$.
Such a net is called an {\it approximate diagonal} for $A$. The
following proposition characterizes approximate
$\varphi$-amenability in terms of approximate diagonal for $A$.

\begin{prop}\label{P:Apd}
 Let $A$ be a Banach algebra and
$\varphi\in\sigma(A)$. Then the following are equivalent:

$(i)$  $A$ is (uniformly) approximately $\varphi$-amenable;

$(ii)$ There exists a net $\{M_\alpha\}\subset (A\widehat{\otimes}
A)^{**}$ such that $(\bigtriangleup^{**}M_\alpha)(\varphi)=1$ and
$||a\cdot M_\alpha-\varphi(a)M_\alpha||\rightarrow 0$ for all $a\in
A$ (uniformly on the unit ball of $A$, respectively);

$(iii)$ There exists a net $\{M_\alpha\}\subset (A\widehat{\otimes}
A)^{**}$ such that
$(\bigtriangleup^{**}M_\alpha)(\varphi)\rightarrow 1$ and $||a\cdot
M_\alpha-\varphi(a)M_\alpha||\rightarrow 0$ for all $a\in A$
(uniformly on the unit ball of $A$, respectively).
\end{prop}

\proofbegin
 $(i)\Rightarrow(ii)$ Choose $u_0\in A$ such that
$\varphi(u_0)=1$. Then $A=\mathbb{C}u_0\oplus I_\varphi$, where
$I_\varphi=\ker\varphi$. We define an action of $A$ on
$X=A\widehat{\otimes}A$ by
$$a\cdot (b\otimes c)=ab\otimes c,  (b\otimes c)\cdot a=\varphi(a)
b\otimes c, ~~a,b,c\in A.$$ Then the dual $A$-bimodule
$(A\widehat{\otimes}A)^*$ is a $(\varphi,A)$-bimodule.

We know that $\varphi\otimes\varphi\in (A\otimes A)^*$, and $a\cdot
(\varphi\otimes\varphi)=(\varphi\otimes\varphi)\cdot
a=\varphi(a)\varphi\otimes\varphi$. Therefore
$\mathbb{C}\varphi\otimes\varphi$ is a closed submodule of
$(A\otimes A)^*$ and $(A\otimes
A)^*/\mathbb{C}\varphi\otimes\varphi$ is a $(\varphi, A)$-bimodule
for which the module actions are given by $a\cdot [f]=\varphi(a)[f],
[f]\cdot a=[f*a]$, for all $a\in A, [f]\in (A\otimes
A)^*/\mathbb{C}\varphi\otimes\varphi$.

 Moreover, define a
derivation $D: A\rightarrow (A\widehat{\otimes}A)^{**} $ given by
$D(a)=au_0\otimes u_0-\varphi(a)u_0\otimes u_0$. Then $D(a)\in
\ker(\varphi\otimes\varphi)\subset
\{\mathbb{C}\varphi\otimes\varphi\}^\bot\cong [(A\otimes
A)^*/\mathbb{C}\varphi\otimes\varphi]^*$. Since $A$ is approximately
$\varphi$-amenable, it follows that there exists a net
$\{N_\alpha\}\subset \{\mathbb{C}\varphi\otimes\varphi\}^\bot$ such
that $D(a)=\lim\limits_\alpha a\cdot N_\alpha -\varphi(a)N_\alpha$
for all $a\in A$.

Set $M_\alpha=u_0\otimes u_0-N_\alpha$, then
$(\bigtriangleup^{**}M_\alpha)(\varphi)=1$ and $||a\cdot
M_\alpha-\varphi(a)M_\alpha||\rightarrow 0$ for all $a\in A$.

$(ii) \Rightarrow(iii)$ It is clear.

$(iii)\Rightarrow(i)$  Now suppose that $(iii)$ holds. Let $E$ be a
$(\varphi, A)$-bimodule and let $D: A\rightarrow E^*$ be a
derivation. For each $\alpha$, set $f_\alpha(x)=M_\alpha(\mu_x)$,
where for $a,b\in A, x\in E, \mu_x(a\otimes b)=(D(a)\cdot b)(x)$.
Then, with $(m_\alpha^\nu)\subset A\widehat{\otimes}A$  converging
$w^{*}$ to $M_\alpha$, and noting that for $m\in
A\widehat{\otimes}A$,
$$\mu_{x\cdot a}(m)=(\mu_x\cdot
a)(m)-(D(a)\cdot\bigtriangleup(m))(x),$$
 we have
\begin{eqnarray*}
(a\cdot f_\alpha)(x)&=&f_\alpha(x\cdot a)\\
&=&M_\alpha(\mu_{x\cdot a})\\
&=& \lim\limits_\nu\mu_{x\cdot
a}(m_\alpha^\nu)\\
&=&M_\alpha(\mu_{x}\cdot
a)-\lim\limits_\nu[D(a)\cdot\bigtriangleup(m_\alpha^\nu)](x)\\
&=&(a\cdot M_\alpha)(\mu_x)-\lim\limits_\nu
\varphi(\bigtriangleup(m_\alpha^\nu))D(a)(x)\\
&=&(a\cdot M_\alpha)(\mu_x)-
(\bigtriangleup^{**}M_\alpha)(\varphi)D(a)(x).
\end{eqnarray*}

Thus, for $a\in A$ and $x\in E$,
\begin{eqnarray*}
||D(a)(x)-(f_\alpha\cdot a-a\cdot f_\alpha)(x)||&\leq&
||D(a)(x)-(\bigtriangleup^{**}M_\alpha)(\varphi)D(a)(x)||\\
&+&||(a\cdot M_\alpha)(\mu_x)-(f_\alpha\cdot
a)(x)||\\
&\leq&|1-(\bigtriangleup^{**}M_\alpha)(\varphi)|\cdot||x||\cdot||D(a)||\\
&+&||(a\cdot M_\alpha)-\varphi(a)M_\alpha||\cdot||\mu_x||,
\end{eqnarray*}
whence $D(a)=\lim\limits_\alpha(f_\alpha\cdot a-a\cdot f_\alpha)$
for all $a\in A$. Set $g_\alpha=-f_\alpha$, then
$D(a)=\lim\limits_\alpha(a\cdot g_\alpha-g_\alpha\cdot a)$ for all
$a\in A$. Since  $D$ was arbitrary, it follows that $A$ is
approximately $\varphi$-amenable.

The proof in the case of uniformly approximately $\varphi$-amenable
is similar. \proofend

The following proposition characterizes (uniform) approximate
$0$-amenability in terms of right approximate identity for $A$.

\begin{prop}\label{P:0-amen}
 Banach algebra $A$ is (uniform) approximate $0$-amenability if
and only if $A$ has a (bounded, respectively) right approximate
identity.
\end{prop}

\proofbegin
 Suppose that $E$ is a $(\varphi, A)$-bimodule  and $D:
A\rightarrow E^*$ is a derivation. Then $D(b)\cdot a=0$, for all
 $a,b\in A$. If $\{e_\alpha\}$ is  (bounded) right approximate identity
for $A$, then $D(a)=\lim\limits_\alpha
D(ae_\alpha)=\lim\limits_\alpha a\cdot D(e_\alpha).$ This implies
that $A$ is (uniformly) approximately $0$-amenable.

The converse is clear from \cite[Lemma 2.1,Lemma 2.2, Theorem
4.2]{G2004}. \proofend

Note that (uniform) approximate right character amenability of $A$
is equivalent to (uniform, respective) approximate
$\varphi$-amenability for all $\varphi\in\sigma(A)$ together with
(uniform, respectively) approximate $0$-amenability. Any statement
about (uniform) approximate right character amenability turns into
an analogous statement about (uniform, respectively) approximate
left character amenability by simply replacing $A$ by its opposite
algebra. Then standard arguments of Proposition \ref{P:P-A1},
\ref{P:RAIA}, \ref{P:Apd} and \ref{P:0-amen} apply, we have the
following theorems.

\begin{thm}\label{T:Ap-a1}
For a Banach algebra $A$ the following are equivalent:

$(i)$ $A$ is (uniformly) approximately character amenable;

$(ii)$ $A$ has (bounded, respectively) both left and right
approximate identities, and for any $\varphi\in \sigma(A)$, there
exist nets $\{ m_\alpha\},\{ m_\alpha^{'}\}\subset A^{**}$ such that
$m_\alpha(\varphi)=1, m_\alpha^{'}(\varphi)=1(or~
m_\alpha(\varphi)\rightarrow1, m_\alpha^{'}(\varphi)\rightarrow1 )$
and $||a\cdot m_\alpha-\varphi(a)m_\alpha||+|| m_\alpha^{'}\cdot
a-\varphi(a)m_\alpha^{'}||\rightarrow 0$ for all $a\in A$ (uniformly
on the unit ball of $A$, respectively);

$(iii)$ $A$ has (bounded, respectively) both left and right
approximate identities, and for any $\varphi\in \sigma(A)$, there
exist nets $\{ m_\alpha\},\{ m_\alpha^{'}\}\subset A$, such that
$\varphi(m_\alpha)=1, \varphi(m_\alpha^{'})=1$ (or~
$\varphi(m_\alpha)\rightarrow1, \varphi(m_\alpha^{'})\rightarrow1 )$
and $||am_\alpha-\varphi(a)m_\alpha||+||
m_\alpha^{'}a-\varphi(a)m_\alpha^{'}||\rightarrow 0$ for all $a\in
A$ (uniformly on the unit ball of $A$, respectively);

$(iv)$  $A$ has (bounded, respectively) both left and right
approximate identities, and for any $\varphi\in \sigma(A)$, there
exist nets $\{M_\alpha\}, \{M_\alpha^{'}\}\subset
(A\widehat{\otimes} A)^{**}$ such that
$(\bigtriangleup^{**}M_\alpha)(\varphi)=1$,
$(\bigtriangleup^{**}M_\alpha^{'})(\varphi)=1
(or~(\bigtriangleup^{**}M_\alpha)(\varphi)\rightarrow 1,
(\bigtriangleup^{**}M_\alpha^{'})(\varphi)\rightarrow 1 )$ and
$||a\cdot M_\alpha-\varphi(a)M_\alpha||+|| M_\alpha^{'}\cdot
a-\varphi(a)M_\alpha^{'}||\rightarrow 0$ for all $a\in A$ (uniformly
on the unit ball of $A$, respectively).
\end{thm}

\begin{thm}\label{T:Ap-a2}
A Banach algebra $A$ is left (right) approximately character
amenable if and only if $\ker\varphi$ has a left (right,
respectively) approximate identity for every
$\varphi\in\sigma(A)\cup\{0\}$.
\end{thm}

\vskip1cm
\section{Hereditary properties of approximate
character amenability}

 In this section, we are concerned with  hereditary
properties of (uniform) approximate character amenability.

\begin{prop}
 Suppose that $A$ is (uniformly) approximately
character amenable (contractible) and $\Phi: A\rightarrow B$ is a
continuous epimorphism. Then $B$ is (uniformly, respectively)
approximately character amenable (contractible, respectively).
\end{prop}

\proofbegin
 The standard argument, \cite[Proposition 2.2]{G2004}
applies. \proofend

\begin{prop}
 Suppose that $A$ is (uniformly) approximately
character amenable (contractible), and $J$ is a closed two-sided
ideal of $A$. Then $A/J$ is (uniformly, respectively) approximately
character amenable (contractible, respectively). If $J$ is character
amenable (contractible) and $A/J$ is (uniformly) approximately
character amenable (contractible), then $A$ is (uniformly,
respectively)approximately character amenable (contractible,
respectively).
\end{prop}

\proofbegin The standard argument, \cite[Proposition 1.30]{A1988}
applies. \proofend

\begin{prop}\label{P:Ideal}
Let $A$ be a Banach algebra and $J$ a weakly complemented ideal of
$A$. Let $\varphi\in \sigma(A)$ satisfy  $\varphi|_{J}\neq 0$. If
$A$ is (uniformly) approximately $\varphi$-amenable, then J is
(uniformly, respectively) approximately $\varphi|_{J}$-amenable.
\end{prop}

\proofbegin Since $A$ is approximately $\varphi$-amenable, it
follows from Proposition \ref{P:P-A1} that there exists a net $\{
m_\alpha\}\subset A^{**}$, such that $m_\alpha(\varphi)=1$ and
$||a\cdot m_\alpha-\varphi(a)m_\alpha||\rightarrow 0$ for all $a\in
A$.

Since $J$ weakly complemented of $A$,  there exist a closed
subspace $X$ of $A^*$ such that $A^*=J^{\bot}\oplus X$. That is to
say, there exists $K>0$ such that for any $F\in A^*$, $F=x_F+y_F$,
where $x_F\in J^{\bot}, y_F\in X$, and  $||x_F||\leq K||F||,
||y_F||\leq K||F||$. If, in addition, $a\in J$,  then $x_F\cdot
a=0$. Thus $||\varphi(a)m_\alpha(x_F)||\rightarrow 0$ for all $a\in
J,$ and uniformly for $F\in A^*$ and $||F||\leq 1$. Choose $a=u_0\in
J$ with $\varphi(u_0)=1$, then $||m_\alpha(x_F)||\rightarrow 0$
uniformly for $||F||\leq 1$.

Set $n_\alpha(f)=m_\alpha(y_F)$ for $f\in J^*$, where $F$ is any
extension of $f$.  Notice  that $||a\cdot
n_\alpha-\varphi(a)n_\alpha||\rightarrow 0, $ and
$n_\alpha(\varphi|_{J})\rightarrow 1$ for all $a\in J$. To see this,
for $a, b\in J, f\in J^*$ and $F$ is a extension of $f$,
$$(y_F\cdot a)(b)=y_F(ab)=f(ab),~ y_{F\cdot a}(b)=f(ab).$$

It follows that there is a $x\in J^\bot$ such that $y_F\cdot
a=x+y_{F\cdot a}$ and $||y_{F\cdot a}||\leq K||y_F\cdot a||,
||x||\leq K||y_F\cdot a||$. Then, for any $f\in J^* $ and $F$ an
extension of $f$,
\begin{eqnarray*}
|a\cdot n_\alpha(f)-\varphi(a)n_\alpha(f)|&=&| n_\alpha(f\cdot
a)-\varphi(a)m_\alpha(y_F)|\\
&=&| m_\alpha(y_{F\cdot a})-\varphi(a)m_\alpha(y_F)|\\
&\leq&|m_\alpha(y_F\cdot a)-\varphi(a)m_\alpha(y_F)|+|
m_\alpha(x)|\\
&\leq&||a\cdot
m_\alpha-\varphi(a)m_\alpha||\cdot||y_F||+|m_\alpha(x)|.
\end{eqnarray*}

 It follows that $||a\cdot
n_\alpha-\varphi(a)n_\alpha||\rightarrow 0$ and
$n_\alpha(\varphi|_{J})=m_\alpha(\varphi)-m_\alpha(x_\varphi)\rightarrow
1$ for all $a\in J$. Then,  by Proposition \ref{P:P-A1}, $J$ is
approximately $\varphi|_{J}$-amenable.

The proof in the case of uniform approximate $\varphi$-amenability
is similar. \proofend

\begin{lem}\label{L: Ext}
Let $A$ be a Banach algebra and $J$ is an ideal of $A$, which has a
right or left approximate  identity. Then every $\varphi\in
\sigma(J)$ can be extended to a $\widetilde{\varphi}$ in
$\sigma(A)$.
\end{lem}

\proofbegin
 Assume that $\{e_\alpha\}$ is right approximate  identity
for $J$ and $\varphi\in \sigma(J)$. Choose $u_0\in J$ such that
$\varphi(u_0)=1$, then $J=\mathbb{C}u_0\oplus I_{\varphi}$, and $
u_0a-u_0\in I_{\varphi}$ for all $a\in J$, where
$I_{\varphi}=\ker\varphi$.

Set $\widetilde{\varphi}(a)=\varphi(u_0a)$ for all $a\in A$. Then
$\widetilde{\varphi}|_{J}=\varphi$, and
$\widetilde{\varphi}(a_1a_2)=\varphi(u_0a_1a_2)=\widetilde{\varphi}(a_1)\widetilde{\varphi}(a_2)$,
since
$u_0a_1a_2-u_0a_1u_0a_2=(u_0a_1-u_0a_1u_0)a_2=\lim\limits_\alpha
(u_0a_1-u_0a_1u_0)e_\alpha a_2\in I_{\varphi}$. It follows that
$\widetilde{\varphi}\in\sigma(A)$.

The proof in the case of a left approximate identity is similar.
\proofend

\begin{thm}\label{T:Ideal}
 Let $A$ be a Banach algebra and $J$ a weakly
complemented ideal of $A$ with  (bounded) left and right approximate
identities. Suppose that $A$ is (uniformly) approximately character
amenable. Then $J$ is (uniformly, respectively)approximately
character amenable.
\end{thm}

\proofbegin Clearly $J$ is (uniformly) approximately
$\varphi$-amenable, for any $\varphi\in\sigma(J)\cup\{0\}$,
 by Proposition \ref{P:0-amen}, \ref{P:Ideal} and Lemma \ref{L: Ext}. That is to say, $J$ is
(uniformly) approximately right character amenable. The proof in the
case of (uniform) approximate left character amenability is similar.
Thus $J$ is (uniformly) approximately character amenable. \proofend

By \cite[Theorem 2.6]{M2008}  a Banach algebra $A$ without a unit is
character amenable if and only if $A^\sharp$ is character amenable.
We obtain a similar result for approximate character amenability.

\begin{lem}\label{L:3}
Let $A$ be a Banach algebra without identity and let $A^\sharp$
denote the unitization of $A$ by adjoining an identity $e$. Let
$\varphi\in \sigma(A)\cup\{0\}$ and let $\varphi_e$ be the unique
extension of $\varphi$ to an element of $\sigma(A^\sharp)$. Then
$A^\sharp$ is approximately $\varphi_e$-amenable if $A$ is
 approximately $\varphi$-amenable.
\end{lem}

\proofbegin Assume that $A$ is approximately $\varphi$-amenable,
where $\varphi\in\sigma(A)$. Then the standard argument of
\cite[Lemma 3.2]{E2008} applies, so that $A^\sharp$ is approximately
$\varphi_e$-amenable.

If $\varphi=0$, $A^\sharp=\mathbb{C}e\oplus A$ and $\varphi_e(A)=0$,
$\varphi_e(e)=1$. Since $A$ is approximately $0$-amenable, it
follows from Proposition \ref{P:0-amen} that $A$ has a right
approximate identity $(b_\alpha)$.

Set $m_\alpha=e-b_\alpha$, then
$||am_\alpha-\varphi_e(a)m_\alpha||\rightarrow 0$ and
$\varphi_e(m_\alpha)=1$, for all $a\in A^\sharp$. Thus $A^\sharp$ is
approximately $\varphi_e$-amenable, by Proposition
\ref{P:P-A1}/\ref{P:Equi 1}. \proofend

\begin{thm}\label{T:unity}
Let $A$ be a Banach algebra without a unit. Then $A$ is
approximately character amenable if and only if $A^\sharp$ is
 approximately character amenable.
\end{thm}

\proofbegin
 Assume $A$ is approximately character amenable, the
 argument of Lemma \ref{L:3} applies,  $A^\sharp$ is
approximately character amenable.

For the converse, assume that $A^\sharp$ is approximately character
amenable, then $A^\sharp$ is both  approximately right and left
character amenable. We shall show that $A$ has both right and left
approximately identities.

Define $\varphi\in\sigma(A^\sharp)$ by  $\varphi(A)=0,
\varphi(e)=1$. Since $A^\sharp$ is right approximately character
amenable, Proposition \ref{P:P-A1}/\ref{P:Equi 1} gives a net
$(m_\alpha)$ in $A^\sharp$ such that $\varphi(m_\alpha)\rightarrow
1$ and $||bm_\alpha||\rightarrow 0$ for all $b\in A $. Writing
$m_\alpha=\lambda_\alpha e+b_\alpha$ where
$\lambda_\alpha\in\mathbb{C}, b_\alpha\in A$, it follows that
$||bb_\alpha+b||\rightarrow 0$ and hence $A$ has a right
approximately identity. Similarly,  $A$ has a left approximately
identity. Then, by Theorem \ref{T:Ideal}, $A$ is approximately
character amenable. \proofend

\vskip1cm
\section{Characterization of (approximate) character contractibility}

In this section, we first characterize $\varphi$-contractibility,
 approximate $\varphi$-contractibility and uniform
approximate $\varphi$-contractibility in two different ways and then
characterize character contractibility,  approximate character
contractibility and uniform approximate character contractibility in
the same different ways.

\begin{prop}\label{P:Phi-C}
Let $A$ be a Banach algebra and $\varphi\in \sigma(A)$. Then the
following conditions are equivalent:

$(i)$ $A$ is $\varphi$-contractible;

$(ii)$ There exists  $m\in A$ such that $\varphi(m)=1$ and
$am=\varphi(a)m$ for all $a\in A$.
\end{prop}

\proofbegin  $(i)\Rightarrow(ii)$
Choose  $u_0\in A$ such that
$\varphi(u_0)=1$. Then $A=\mathbb{C}u_0\oplus I_\varphi$, where
$I_\varphi=\ker\varphi$. We define an action of $A$ on $X=A$ by
$a\cdot x=ax, x\cdot a=\varphi(a)x, a\in A, x\in X$. Then $X$ is a
$(A, \varphi)$-bimodule.

Moreover, define a derivation $D: A\rightarrow A$ given by
$D(a)=au_0-\varphi(a)u_0$. Then $D(a)\in \ker\varphi$ and
$\ker\varphi$ is a submodule of $X$. Since $A$ is
$\varphi$-contractible, it follows that there exists $n\in
\ker\varphi$ such that $D(a)=a\cdot n-n\cdot a$. Set $m=u_0-n$. Then
$\varphi(m)=1$ and $am=\varphi(a)m$ for all $a\in A$ .

$(ii)\Rightarrow(i)$   $m\in A$ is such that $\varphi(m)=1$ and
$am=\varphi(a)m$ for all $a\in A$. Let $X$ be a $(A,
\varphi)$-bimodule, and let $D: A\rightarrow X$ be a derivation. Set
$x=D(m)$, then
\begin{eqnarray*}
a\cdot D(m)&=&D(am)-D(a)\cdot m\\
&=&D(am)- \varphi(m)D(a)\\
&=&\varphi(a)D(m)-D(a).
\end{eqnarray*}

It follows that $D(a)=\varphi(a)D(m)-a\cdot D(m)=D(m)\cdot a-a\cdot
D(m)$ for all $a\in A$. Since  $D$ was arbitrary, it follows that
$A$ is
 $\varphi$-contractible.
\proofend

\begin{prop}\label{P:P-C1}
 For a Banach algebra $A$ the following are
equivalent:

$(i)$ $A$ is approximately $\varphi$-contractible;

(ii) There exists  a net $\{ m_\alpha\}\subset A$ such that
$\varphi(m_\alpha)=1$ and
$||am_\alpha-\varphi(a)m_\alpha||\rightarrow 0$ for all $a\in A$;

$(iii)$ There exists  a net $\{ m_\alpha\}\subset A$ such that
$\varphi(m_\alpha)\rightarrow1$ and
$||am_\alpha-\varphi(a)m_\alpha||\rightarrow 0$ for all $a\in A$.
\end{prop}

\proofbegin The proof is a minor modification of the proof of the
analogous statement in Proposition \ref{P:Phi-C}. \proofend

\begin{lem}\label{L:C1}
 Let $A$ be a Banach algebra and
$\varphi\in\sigma(A)$. If the ideal $I_\varphi=\ker\varphi$ has a
right  identity, then $A$ is $\varphi$-contractible.
\end{lem}

\proofbegin Choose  $u_0\in A$ such that $\varphi(u_0)=1$. Then
$A=\mathbb{C}u_0\oplus I_\varphi$ and there exists $a_0\in
I_\varphi$ such that $u_0^2-u_0=a_0$. Let $b_0$ be a right identity
for $I_\varphi$. Set $m=u_0-u_0b_0\in A$. Then, for any $b\in
I_\varphi$, $b(u_0-u_0b_0)=bu_0-bu_0b_0= 0.$ Furthermore,
$$u_0(u_0-u_0b_0)-(u_0-u_0b_0)$$
$$=u_0^2-u_0^2b_0-u_0+u_0b_0$$
$$=a_0-a_0b_0= 0.$$
It follows that $am-\varphi(a)m= 0$ for all $a\in A$ and
$\varphi(m)=1$. Thus $A$ is $\varphi$-contractible by Proposition
\ref{P:Phi-C}. \proofend

\begin{lem}\label{L:C2}
 Suppose that $A$ is
$\varphi$-contractible for some $\varphi\in\sigma(A)$ and that $A$
has a right identity. Then $I_\varphi=\ker \varphi$ has a right
identity.
\end{lem}

\proofbegin Assume that $A=\mathbb{C}u_0\oplus I_\varphi$ for some
$u_0\in A$, and $n=u_0+b_0$ be a right identity for $A$, where
$b_0\in I_\varphi$. Since $A$ is  $\varphi$-contractible, it follows
from Proposition \ref{P:Phi-C} that there exist an $m=u_0+b_1\in A$
such that $am-\varphi(a)m= 0$, where $b_1\in I_\varphi$. Set
$e=b_0-b_1$. Then $be-b= 0,$ for any $b\in I_\varphi$. In fact, for
any $b\in I_\varphi$, $bn-b= 0$, and $bm= 0$. It follows that
$bu_0+bb_0-b= 0$ and $bu_0+bb_1= 0$. Thus $be-b=0$ for all $b\in
I_\varphi$. \proofend

The following result follows immediately from Lemma
\ref{L:C1}/\ref{L:C2} and we omit its proof.
\begin{prop}\label{P:RAI}
 Let $A$ be a Banach algebra with a right
 identity and $\varphi\in\sigma(A)$. Then $A$ is
$\varphi$-contractible if and only if $I_\varphi$ has a right
identity.
\end{prop}

For $0$-contractibility and (uniform) approximate
$0$-contractibility we have the following parallel results to
approximate $0$-amenability, the proofs of Proposition \ref{P:0-con}
and Proposition \ref{P:A0-con} are minor modifications of the proof
of the analogous statements in Proposition \ref{P:0-amen} and will
be omitted.

\begin{prop}\label{P:0-con}
Banach algebra $A$ is $0$-contractible  if and only if $A$ has a
right identity.
\end{prop}

\begin{prop}\label{P:A0-con}
Banach algebra $A$ is approximately $0$-contractible if and only if
$A$ has a right approximate identity.
\end{prop}

\begin{prop}\label{P:UA0-con}
Banach algebra $A$ is uniformly approximately $0$-contractible  if
and only if $A$ has a right identity.
\end{prop}

\proofbegin It suffices to show that if $A$ is uniformly
approximately $0$-contractible, then $A$ has a right identity.
Define an action of $A$ on $E=A$ by $a*x=ax, x* a=0,~ a\in A, x\in
E,$ then $E$ is $(A, 0)$-bimodule. The natural injection $a\mapsto a
: A\rightarrow A$ is a derivation. Thus there is a net
$\{e_\alpha\}$ in $A$ such that $ae_\alpha\rightarrow a$ uniformly
for $||a||\leq 1$.  Let $R_b$ denote right multiplication by $b\in
A$. Then there is $e_\lambda\in \{e_\alpha\}$  with
$||R_{e_\lambda}a-a||<||a||$ for all $a\in A$. Thus $R_{e_\lambda}$
is invertible. It follows that, there is $e\in A$ such that
$ee_\lambda=e_\lambda$, whence $(ae-a)e_\lambda=0$ for all $a\in A$.
Then $e$ is a right identity of $A$ by injectivity of
$R_{e_\lambda}$. \proofend

Finally, we characterize character contractibility and (uniform)
approximate character contractibility, whose proofs are minor
modifications of the proofs of the analogous statements in Theorem
\ref{T:Ap-a1}/\ref{T:Ap-a2} and will be omitted.

\begin{thm}
 For a Banach algebra $A$ the following are equivalent:

$(i)$ $A$ is character contractible;

 $(ii)$ $A$ has an identity and for any $\varphi\in \sigma(A)$, there
exists $ m_1,  m_2\in A$ such that $m_i(\varphi)=1, (i=1, 2)$ and
$am_1-\varphi(a)m=0$, $m_2a-\varphi(a)m_2=0$ for all $a\in A$;

$(iii)$ $A$ has an identity and for any $\varphi\in \sigma(A)$,
there exists bounded nets $\{M_\alpha\}, \{M_\alpha^{'}\}\subset
(A\widehat{\otimes} A)$ such that $\varphi(\bigtriangleup
M_\alpha)=1$, $\varphi(\bigtriangleup M_\alpha^{'})=1
(or~\varphi(\bigtriangleup M_\alpha)\rightarrow 1,
\varphi(\bigtriangleup M_\alpha^{'})\rightarrow 1 )$ and
 $||a\cdot
M_\alpha-\varphi(a)M_\alpha||+|| M_\alpha^{'}\cdot
a-\varphi(a)M_\alpha^{'}||\rightarrow 0$ for all $a\in A$.
\end{thm}

\begin{thm}\label{T:con 2}
Let $A$ be a Banach algebra. Then $A$ is  character contractible if
and only if $\ker\varphi$ has an identity for every
$\varphi\in\sigma(A)\cup\{0\}$.
\end{thm}

\begin{thm}
 For a Banach algebra $A$ the following are equivalent:

$(i)$ $A$ is (uniformly)  approximately character contractible;

$(ii)$ $A$ has (an identity, respectively) both right and left
approximate identities  and for any $\varphi\in \sigma(A)$, there
exist nets $\{ m_\alpha\},\{ m_\alpha^{'}\}\subset A$ such that
$m_\alpha(\varphi)=1, m_\alpha^{'}(\varphi)=1(or~
m_\alpha(\varphi)\rightarrow1, m_\alpha^{'}(\varphi)\rightarrow1 )$
and $||am_\alpha-\varphi(a)m_\alpha||+||
m_\alpha^{'}a-\varphi(a)m_\alpha^{'}||\rightarrow 0$ for all $a\in
A$ (uniformly on the unit ball of $A$, respectively);

$(iii)$ $A$ has (an identity, respectively) both right and left
approximate identities  and for any $\varphi\in \sigma(A)$, there
exist nets $\{M_\alpha\}, \{M_\alpha^{'}\}\subset
(A\widehat{\otimes} A)$ such that $\varphi(\bigtriangleup
M_\alpha)=1$, $\varphi(\bigtriangleup M_\alpha^{'})=1
(or~\varphi(\bigtriangleup M_\alpha)\rightarrow 1,
\varphi(\bigtriangleup M_\alpha^{'})\rightarrow 1 )$ and
 $||a\cdot
M_\alpha-\varphi(a)M_\alpha||+|| M_\alpha^{'}\cdot
a-\varphi(a)M_\alpha^{'}||\rightarrow 0$ for all $a\in A$ (uniformly
on the unit ball of $A$, respectively).
\end{thm}


\vskip1cm
\section{Relations between generalized notions of character amenability}

  In this section, we are concerned with relations among generalized
notions of character amenability (contractibility). Firstly, we
prove that $w^*$-approximate $\varphi$-amenability, approximate
$\varphi$-amenability and approximate $\varphi$-contractibility are
the same properties. We then prove that  approximate character
amenability and approximate character contractibility are the same
properties, as are uniform approximate character amenability
(contractibility) and character amenability (contractibility,
respectively). Moreover, we obtain that character contractibility
and contractibility are the same properties for commutative Banach
algebra.

\begin{lem}\label{L:Equi 2}
 For a Banach algebra $A$ and
$\varphi\in\sigma(A)$, the following are equivalent:

$(i)$ A is approximately $\varphi$-contractible;

$(ii)$ A is approximately $\varphi$-amenable;

$(iii)$ A is $w^*$-approximately $\varphi$-amenable.
\end{lem}

\proofbegin It suffices to show that $(iii)\Rightarrow (i)$. Suppose
that $(iii)$ holds. By the standard argument of Proposition
\ref{P:P-A1}, it follows that there is a net $\{m_\alpha\}\subset
A^{**}$ such that $(a\cdot
m_\alpha-\varphi(a)m_\alpha)(f)\rightarrow 0,
m_\alpha(\varphi)\rightarrow 1$ for all $a\in A, f\in A^*$. It
follows from the proof of
 Proposition \ref{P:Equi 1} that there exists a net $\{ n_\beta\}\subset A$
such that $ \varphi(n_\beta)\rightarrow1$ and $||an_\beta-\varphi(a)
n_\beta||\rightarrow 0$ for all $a\in A$.
 Thus A is approximately $\varphi$-contractible by Proposition \ref{P:P-C1}.
 \proofend

\begin{thm}\label{T:Equi 1}
A Banach algebra $A$ is approximately character contractible if and
only if  $A$ is approximately character amenable.
\end{thm}

\proofbegin
 The standard argument of Proposition \ref{P:0-amen}/\ref{P:A0-con} and
Lemma \ref{L:Equi 2} applies. \proofend

Recall that a Banach algebra $A$ is uniformly approximately
$\varphi$-amenable, if  for every  $(\varphi, A)$-bimodule $E$,
every derivation $D$ from $A$ into the dual $A$-bimodule $E^*$ may
be approximated uniformly on the unit ball of A by inner
derivations. Clearly any $\varphi$-amenable ($\varphi$-contractible)
Banach algebra is uniformly approximately $\varphi$-amenable
($\varphi$-contractible, respectively). In the following  theorems
we prove that the converse is also true.  That is to say, uniform
approximate character amenability is equivalent to character
amenability, and uniform approximate character contractibility is
equivalent to character contractibility.

\begin{lem}\label{L:4}
\cite[Corollary 2.7]{M2008}  Let $A$ be a Banach algebra, if
$\ker\varphi$ has both bounded left and right approximate identities
for every $\varphi\in \sigma(A)\cup\{0\}$. Then $A$ is character
amenable.
\end{lem}

\begin{prop}\label{P:Uap-a}
 Let $A$ be unital Banach algebra, then $A$
is uniformly approximately character amenable if and only if it is
character amenable.
\end{prop}

\proofbegin It suffices to show that if $A$ is uniformly
approximately character amenable then $A$ is character amenable.

Take $\varphi\in \sigma(A)$. Since $A$ is uniformly approximately
right character amenable it follows from Proposition \ref{P:P-A1}
that there is a net $\{m_\alpha\}\subset A^{**}$, such that
$m_\alpha(\varphi)=1$ and $||a\cdot
m_\alpha-\varphi(a)m_\alpha||\rightarrow 0$ uniformly for $
||a||\leq 1$. Note that $A=\mathbb{C}e\oplus I_{\varphi}$, where
$I_{\varphi}=\ker\varphi$, and let $J(\varphi)= \{n \in A^{**},
n(\varphi)= 0\}$.  $J(\varphi)$ is a $w^*$-closed ideal of $A^{**}$,
and $J(\varphi)$ can be canonically identified with the second dual
$(I_\varphi)^{**}$. We have $||a\cdot m_\alpha||\rightarrow
 0$ uniformly on the unit ball of $I_\varphi$.

Set $n_\alpha=e-m_\alpha$. Then $n_\alpha(\varphi)=0$, that is,
$n_\alpha\in J(\varphi)$. Moreover, for any $a\in I_{\varphi}$,
$|(an_\alpha-a)(f)|=|(a\cdot n_\alpha-a)(f)|=|(a\cdot
m_\alpha)(f)|\leq||a\cdot m_\alpha||\cdot||f||$. It follows that
$||an_\alpha-a||\rightarrow 0$ uniformly on the unit ball of
$I_\varphi$. Now take $s\in (I_\varphi)^{**}$,  then,  there is a
net $(s_i )\subset I_\varphi$ such that $||s_i||\leq ||s|| $ and
$s_i\rightarrow s(w^*~in ~i)$. Thus $s_in_\alpha-s_i\rightarrow
sn_\alpha-s(w^* ~in ~i)$ and $||sn_\alpha-s||\leq
\sup_i||s_in_\alpha-s_i||$. It follows that
$||an_\alpha-a||\rightarrow 0$ uniformly for $ a\in
(I_{\varphi})^{**}$ and $||a||\leq 1$.

Thus there is a sequence $(n_k)\subset(I_\varphi)^{**}$  and
$\varepsilon_k\rightarrow 0$ such that
$$||an_k-a||\leq \varepsilon_k||a||, (a\in(I_\varphi)^{**}).$$

Thus, the multiplication operator $R_{n_k}:
(I_\varphi)^{**}\rightarrow (I_\varphi)^{**}$ defined by
$R_{n_k}(s)=sn_k$ satisfies $||R_{n_k}-id_{(I_\varphi)^{**}}||< 1$
for $k$ sufficiently large. Take such $k$, so that $R_{n_k}$ is
invertible. By surjectivity, there is $x\in(I_\varphi)^{**}$ such
that $xn_k = n_k$. Then for each $y\in(I_\varphi)^{**}$ we have
$(yx-y)n_k = 0$. From the injectivity of $R_{n_k}$ this implies $yx
= y$ for all $y\in(I_\varphi)^{**}$. So $(I_\varphi)^{**}$ has a
right identity, then $I_\varphi$ has a bounded right approximately
identity.

 On the other hand, since $A$ is uniformly approximately left  character
 amenable.   We can deduce that
 $I_\varphi$ has a bounded left approximately identity. Since
 $\varphi$ is arbitrary, it follows from Lemma \ref{L:4} that
 $A$ is character amenable.
\proofend

\begin{cor}\label{C:1}
If $A^\sharp$ is uniformly approximately character amenable, then
$A$ has bounded right and left approximate identity.
\end{cor}

\proofbegin Choose $\varphi\in \sigma(A^\sharp)$ such that
$\varphi(A)=0, \varphi(e)=1$, the standard argument of Proposition
\ref{P:Uap-a} applies. \proofend

\begin{lem}\label{L:5}
  Let $A$ be a Banach algebra without a unit. Then
$A$ is uniformly approximately character amenable if and only if
$A^\sharp$ is uniformly approximately character amenable.
\end{lem}

\proofbegin
 Assume $A$ is uniformly approximately character
amenable, the standard argument of Lemma \ref{L:3} applies,
$A^\sharp$ is uniformly approximately character amenable.

For the inverse, assume $A^\sharp$ is uniformly approximately
character amenable. It follows from Corollary \ref{C:1} that $A$ has
bounded left and right approximately identity. Thus $A$ is uniformly
approximately character amenable by Theorem \ref{T:Ideal}. \proofend

Note that $A$ is (uniformly approximately) character amenable if and
only if its unitization $A^\sharp$ is (uniformly approximately,
respectively) character amenable  \cite[Theorem 2.6]{M2008} and
Lemma \ref{L:5}, then, by Proposition \ref{P:Uap-a}, we have the
following theorem.

\begin{thm}\label{T: Uaa-a}
  A Banach algebra $A$ is uniformly approximately
character amenable if and only if it is character amenable.
\end{thm}

\begin{cor} If a finite-dimensional Banach algebra is uniformly
approximately character amenable, then it is  character amenable.
\end{cor}

For uniform approximate character contractibility and character
contractibility , we have the following parallel result  whose proof
is a minor modification of the proof of the analogous statements in
Proposition \ref{P:Uap-a} and will be omitted.

\begin{thm}\label{T: Uac-c}
 A Banach algebra $A$ is uniformly approximately
character contractible if and only if it is character contractible.
\end{thm}

\begin{cor}
If a finite-dimensional Banach algebra is approximately character
contractible, then it is  character contractible.
\end{cor}

Now, we shall conclude this section by proving the following result.

\begin{thm}\label{T: c-cc}
 Let $A$ be a commutative Banach algebra. Then $A$ is
character contractible if and only if $A$ is  isomorphic to
$\mathbb{C}^n$.
\end{thm}

\proofbegin It suffices to show that if $A$ is character
contractible then $A$ is  contractible.

 Let $\mathfrak{M}$ be the maximal ideal space of $A$. It follows from Theorem \ref{T:con 2}
 that, given $\varphi\in\mathfrak{M}\cup\{0\}$, there exists $E_\varphi\in
\ker\varphi$ such  that $E_\varphi$ is an identity for
$\ker\varphi$. Arbitrarily choose $\varphi_{1},
\varphi_{2}\in\mathfrak{M}$. If $\varphi_{1}\neq\varphi_{2}$, then
$\varphi_{1}(E\varphi_2)=1$ and $\varphi_{2}(E\varphi_1)=1$. It
follows that each point of $\mathfrak{M}$ is isolated, so that
$\mathfrak{M}$ is finite since $\mathfrak{M}$ is compact. Hence
$\mathfrak{M}=\{\varphi_1, \varphi_2, \cdots, \varphi_n\}$ and
$E_{\varphi_i}$ is the identity for $\ker\varphi_i$ for all $1\leq
i\leq n$. Hence, $E_0=E_{\varphi_1}\cdot E_{\varphi_2}\cdot
\cdots\cdot E_{\varphi_n}$ is an identity for
Rad$(A)=\ker\varphi_1\cap\cdots\cap \ker\varphi_n$. This is possible
only if Rad$(A)=\{0\}$, so that $A$ is semisimple. It follows that
$A\cong \mathbb{C}^n$, then, by \cite{V1774}, $A$ is contractible.
\proofend


\vskip1cm
\section{Examples}

In this section, we give three examples. The first example shows
that there exists a Banach algebra which is approximately character
amenable but not uniformly approximately character amenable. The
second example shows that there exists  a Banach algebra which is
character amenable but not character contractible. The last example
shows that there exists a Banach algebra which is character
contractible but not amenable.

Given a family $(A_\alpha)_{\alpha\in\Lambda}$ of Banach algebras
defined their $l^{\infty}$ direct sum as
$$l^{\infty}(A_\alpha)=\{(a_\alpha): a_\alpha\in A_\alpha,
||(a_\alpha)||=\sup||a_\alpha||<\infty\}.$$
 Further, set
$$c_0(A_\alpha)=\{(a_\alpha): (a_\alpha)\in l^{\infty}(A_\alpha), ||a_\alpha||\rightarrow
0\}.$$

{\bf Example 1.} There exists an approximately character amenable
Banach algebra which is not uniformly approximately character
amenable.

For each $n\in \mathbb{N}$, take $A_n=\mathbb{C}^n$, each with the
corresponding $l^1$ norm. Then $A_n$ has an identity
$e_n=(1,1,\cdots, 1)$ of norm $n$;
 clearly each $A_n$ is
amenable, so $c_0(A_n^\sharp)$ is approximately amenable by
\cite[Example 6.1]{G2004}, then it is also approximately character
amenable. If $f_\alpha$ is a left approximate identity for
$c_0(A_n^\sharp)$ then, given $n$, there is $\alpha$ such that
$||f_\alpha e_n-e_n||\leq 1$. Hence we have
$||f_\alpha||\geq||f_\alpha e_n||\geq n-1$ and then $\{f_\alpha\}$
is unbounded.

But $c_0(A_n^\sharp)$ is not uniformly approximately character
amenable. Indeed, if it were, Theorem \ref{T:Ap-a1} would imply that
$c_0(A_n^\sharp)$ had  bounded left and right approximate identity.

{\bf Example 2.} There exists a character amenable Banach algebra
which is neither amenable nor
 character contractible.

 Let $V$ be the Volterra operator, $A_V$ be the Banach algebra
 generated by $V$ and $A$ be the Banach algebra
 generated by $V$ and $e$. Then $A=\mathbb{C}e\oplus A_V$, and $A_V$
 has a  bounded two sided approximate identity  \cite[Theorem 5.10]{D1988}. And
 $\sigma(A)=\{\varphi\},$ here $\varphi(A_V)=0, \varphi(e)=1$.
Thus $A$ is a character amenable by Lemma \ref{L:4}. However, $A$ is
not amenable by \cite{D2005}.

Moreover, $A$ is not character contractible. Indeed, if it were,
Proposition \ref{P:Phi-C} implies that there exists an $m=e+b\in A$
such that $a(e+b)=\varphi(a)(e+b)$ for all $a\in A$. If, in
addition, $a\in A_V$ and $a\ne 0$, then $a(e+b)=0$; since $e+b$ is
invertible, we obtain $a=0$, a contradiction.

{\bf Example 3.} Let $A=B(H)$, where $H$ is an infinite-dimensional
Hilbert space. Then $A$ has an identity and $A$ is character
contractible by Theorem \ref{T:con 2}. But $A$ is not amenable by
\cite{V1774}.

\end{document}